\documentclass[14pt]{amsart}
\usepackage{amsfonts}
\usepackage{bm}
\usepackage{amsfonts,amsmath,amssymb,amscd,bbm,amsthm,mathrsfs,dsfont}
\usepackage{mathrsfs}
\usepackage{pb-diagram}
\usepackage{color}
\usepackage{amssymb}
\usepackage{xypic}
\usepackage{amsmath}
\usepackage{amssymb}
\usepackage{amsthm}
\usepackage{amscd}
\usepackage{indentfirst}
\usepackage{graphicx}
\usepackage{epstopdf}
\usepackage{float}
\usepackage{caption}
\usepackage{subcaption}
\usepackage{refstyle}
\usepackage{fancyhdr}
\usepackage[dvipsnames]{xcolor}

\newtheorem{mytheo}{Theorem}[section]

\newtheorem{cons}[mytheo]{Construction}
\newtheorem{cor}[mytheo]{Corollary}

\newtheorem{lem}[mytheo]{Lemma}

\newtheorem{re}[mytheo]{Remark}

\newtheorem{dfn}[mytheo]{Definition}

\newtheorem{eg}[mytheo]{Example}

\title[maximal green sequences for acyclic sign-skew-symmetric matrices]
{sign-coherence of c-vectors and maximal green sequences for acyclic sign-skew-symmetric matrices}
\author{Diana Ahmad$\;\;\;\;\;\;$Fang Li}
\address{Diana Ahmad
\newline Department of Mathematics, Zhejiang University (Yuquan Campus), Hangzhou, Zhejiang 310027, P. R. China}
\email{diana.ahmadsy@yahoo.com}
\address{Fang Li
\newline Department of Mathematics, Zhejiang University (Yuquan Campus), Hangzhou, Zhejiang 310027, P. R. China}
\email{fangli@zju.edu.cn}
\date{version of \today}
\begin{document}

\maketitle

\renewcommand{\thefootnote}{\alph{footnote}}
\setcounter{footnote}{-1} \footnote{\emph{Corresponding author}: Fang Li, E-mail: fangli@zju.edu.cn.}
\renewcommand{\thefootnote}{\alph{footnote}}
\setcounter{footnote}{-1} 
\begin{abstract}
In this paper we construct an unfolding for $c-$vectors of acyclic sign-skew symmetric matrices and we also prove that the sign-coherence property holds for acyclic sign-skew-symmetric matrices. Then we prove that every acyclic sign-skew-symmetric matrix admits a maximal green sequence.
\end{abstract}
\section{introduction and preliminaries}
\numberwithin{equation}{section}
 The problem posed by  A. Berenstein, S. Fomin and A. Zelevinsky in \cite{nw} on whether any acyclic sign-skew-symmetric integer matrix is totally sign-skew-symmetric or not, was a great motivation for many mathematicians to study such matrices. M. Huang and F. Li gave an affirmative answer to this problem and proved in \cite{mingli} that  acyclic sign-skew-symmetric matrices are totally mutable. The authors in \cite{mingli} also proved that every acyclic sign-skew-symmetric matrix can be covered by an (infinite) skew-symmetric matrix which is represented by an (infinite) cluster quiver and this covering can perform arbitrary steps of orbit-mutations. This (infinite) quiver is called an unfolding of this acyclic sign-skew-symmetric matrix. The existence of {such an} unfolding quiver for every {acyclic} sign-skew-symmetric matrix allows us to tackle problems {related to} an {acyclic} sign-skew-symmetric cluster algebra by promoting these problems to an (infinite) skew-symmetric cluster algebra. In this note we {try to find an unfolding} for the $c-$vectores of an acyclic sign-skew-symmetric matrix and prove that it always exists {(see paragraph 3 of  Remark \ref{unfolding2})}. In other words, we prove that every extended acyclic sign-skew-symmetric matrix $\widetilde{B}=\begin{pmatrix}B\\I_n\end{pmatrix}$ can be covered by an (infinite) extended skew-symmetric matrix $\widetilde{B^\S}=\begin{pmatrix}B^\S\\I_\infty\end{pmatrix}$. The construction of this {covering} keeps the principal part as it was constructed in \cite{mingli} which makes the ability of performing arbitrary steps of orbit-mutations {remain} valid. Using the unfolding method, we prove {that the sign-coherence property holds} for acyclic sign-skew-symmetric matrices, which together with {the fact that these matrices are totally mutable}, {means} that maximal green sequences are well-defined for acyclic sign-skew-symmetric matrices. Finally we prove that every acyclic sign-skew-symmetric matrix admits a maximal green sequence (see Theorem \ref{mainn}).\\

A {\bf{skew-symmetric}} matrix is an integer matrix $B=(b_{ij})$ of the size $n \times n$, such that $b_{ij}=-b_{ji}$ for all $1 \leq i,j \leq n$. A {\bf{skew-symmetrizable matrix}} is an integer matrix $B=(b_{ij})$ of the size $n \times n$, such that $B=-(BD)^T$ for $D$ is a diagonal matrix with positive integers. $D$ is called the symmetrizing matrix. A {\bf{sign-skew-symmetric}} matrix is an integer matrix $B=(b_{ij})$ of the size $n \times n$, such that either $b_{ij},b_{ji}= 0$ or $b_{ij}.b_{ji}< 0$ for any $1 \leq i,j \leq n$.\\
{The mutation of a matrix ${B}$ in direction $k$ where $1 \leq k \leq n$  is the matrix}
	 $\mu_k(B)=B^{\prime}=(b_{ij}^{\prime})$ where:\begin{equation}
	\label{firstmutation}
	b_{ij}^{\prime}=\begin{cases} -b_{ij}\,\,\,\,\,,\,\,\,\,\,\,\,\,\,\,\,\,\,\,\,\,\,\,\,\,\,\,\,\,\,\,\,\,\,\,\,\,\,\,\,\,\,\,\,\,\,\,\,\,\,\,\,\,\,\,\,\,\,\,\,\,\,\,\,\,\,\, \text{if}\,\,\, i \,\,\, or \,\,\, j=k \\
	 b_{ij}+\dfrac{1}{2}(\mid{b_{ik}}\mid b_{kj}+b_{ik} \mid b_{kj} \mid)\,\,\,\,\,\text{otherwise}
	 \end{cases}
	 \end{equation}
{Equation (\ref{firstmutation})} is called the {{\bf{matrix mutation formula}}}. The mutation is an involution i.e $\mu_{k}\mu_k(B)=B$.
{ A skew-symmetric matrix $B=(b_{ij}) \in Mat_{n \times n}(\mathbb{Z})$ can be represented by a directed diagram called a {\bf{quiver}} with $n$ vertices such that there are {{$\mid b_{ij} \mid$ many arrows from $j$ to $i$}} if $b_{ij}\geq 0$. $Q_0$ is the set of vertices in $Q$ and $Q_1$ is the set of arrows in $Q$. The mutation formula can be translated to the language of quivers such that
 for every $k \in Q_0$, the {\bf{quiver mutation}} in direction $k$ is obtained by the following steps }
 {\begin{enumerate}
     \item for each subquiver $i\rightarrow k \rightarrow j$ add a new arrow {$i \rightarrow j$.}
     \item  reverse all arrows with source or target $k$.
     \item  remove the arrows in a maximal set of pairwise disjoint 2-cycles.
   \end{enumerate}}

{$Q$ is finite if $Q_0$ and $Q_1$ are both finite. A vertex $i$ falls in the neighbourhood of a vertex $j$ if there is an arrow connecting $i$ and $j$.}

We can easily check that the skew-symmetricity and the skew-symmetrizablity are invariant under mutation, whereas the sign-skew-symmetricity is not necessarily invariant under mutation.
 A sign-skew-symmetric matrix which remains sign-skew-symmetric under any arbitrary finite sequence of mutation is called {\bf{totally sign-skew-symmetric matrix}}. \\
{ An $n\times n$ sign-skew-symmetric matrix $B$ can be associated with a (simple) quiver $\Delta(B)$
with vertices $1, \cdots , n$ such that for {each pair $(i, j)$ with $b_{ij}<0$}, there is exactly one
arrow from vertex $i$ to vertex $j$. Trivially, $\Delta(B)$ has no loops and no 2-cycles.
Recall that the sign-skew-symmetric matrix $B$ is called {\bf{acyclic}} if $\Delta(B)$ is acyclic i.e, $\Delta(B)$ does not admit any directed cycles \cite{mingli}.}
{
\begin{dfn}
\rm{Let $B$ be a totally sign-skew-symmetric matrix, we call $\widetilde{B}=\begin{pmatrix} B\\I_n\end{pmatrix}\in Mat_{2n \times n}(\mathbb{Z})$ the {\bf{extended matrix }}of $B$. And let $\widetilde{B}^{\sigma_m}=\begin{pmatrix} B^{\sigma_m}\\C^{\sigma_m}\end{pmatrix}$ be the matrix obtained from $\widetilde{B}$ by a composition of mutations $\mu_{\sigma_m}=\mu_{k_m}\mu_{k_{m-1}}.....\mu_{k_0}$ such that $1 \leq k_j \leq n$ for $0 \leq j \leq m$. Then the lower part of $\widetilde{B}^{\sigma_m}$ is called the {\bf{$C-$matrix}} and its columns are called the {\bf{$c-$vectors.}}}
\end{dfn}
{{The mutation of a matrix $\widetilde{B}$ in direction $k$ where $1 \leq k \leq n$  is the matrix}
	 $\mu_k(\widetilde{B})=\widetilde{B^{\prime}}=\begin{pmatrix}B^\prime\\C^\prime \end{pmatrix}$ where $B^\prime$ is given as in Equation (\ref{firstmutation}) and $C^\prime=(c_{ij}^{\prime})$ such that:\begin{equation}
	\label{cmutation}
	c_{ij}^{\prime}=\begin{cases} -c_{ij}\,\,\,\,\,,\,\,\,\,\,\,\,\,\,\,\,\,\,\,\,\,\,\,\,\,\,\,\,\,\,\,\,\,\,\,\,\,\,\,\,\,\,\,\,\,\,\,\,\,\,\,\,\,\,\,\,\,\,\,\,\,\,\,\,\,\,\, \text{if} \,\,\, j=k \\
	 c_{ij}+\dfrac{1}{2}(\mid{c_{ik}}\mid b_{kj}+c_{ik} \mid b_{kj} \mid)\,\,\,\,\,\text{otherwise}
	 \end{cases}
	 \end{equation}}
{\begin{re}
\rm{In this paper we refer to the mutation {{given in Equations}} {(\ref{firstmutation}) and (\ref{cmutation})} as {\bf{ordinary mutation}} and the mutation {{given in}} Equation {(\ref{eq1})} as {\bf{orbit-muation}}}.
\end{re}}
By convention $\mu_{k_0}(\widetilde{B})=\widetilde{B}$.({{$\mu_{k_0}$} means no mutation has been applied yet and any $c-$vector in $\widetilde{B}$ has its entries either all non-positive or all non-negative}.)\\
If {the entries} of any $c-$vector in the matrix $\widetilde{B}^{\sigma_j}$ such that $0\leq j < \infty$ are either all non-positive or all non-negative, then we say that the {\bf{sign-coherence property}} {for $C-$matrix} holds for the matrix $\widetilde{B}$. }\\
The idea of the unfolding method of an {acyclic} sign-skew-symmetric matrix $B$ is to create an {(infinite)} quiver $Q$ which covers $B$ and can do orbit-mutations. We recall the way to create such quiver as it was mentioned in \cite{mingli}.\\
A {\bf{locally-finite}} quiver is an infinite quiver which has finitely many arrows incident to each of its vertices. A locally-finite quiver $Q$ can be represented by an {infinite and well-defined matrix $B^\S=(b^\S_{ij})$ {{called the adjacency matrix}} of $Q$} such that $b^\S_{ij}\geq 0$ if there are $\mid b^\S_{ij} \mid$ many arrows from {$j$ to $i$} in $Q$.
{\begin{dfn}
\rm{ Let $B^\S$ be the {adjacency matrix of} a locally-finite quiver $Q$ and let $g$ be a permutation {acting} on $Q_0$, then $g$ is said to be an {\bf automorphism} of $B$ or an {{{\bf{automorphism}}} of $Q$ if
\begin{center}
	{$b_{gi,gj}^\S = b_{ij}^\S \;\;\text{for every} \;\;  i,j \in Q_0$}
\end{center}}}
\end{dfn}
{Let $Q$ be a locally-finite quiver and $\Gamma$ be a subgroup of the symmetric group $S_{Q_0}$}. If all the elements of $\Gamma$ are automorphisms of $Q$, then $\Gamma$ is said to be a {\bf{group of automorphisms}} of this quiver.}
  Let $Q$ be a locally-finite quiver equipped with a group of automorphisms $\Gamma$.{ We denote the orbits created under the action of $\Gamma$ by $\bar{i}$ such that {$i \in Q_0$}. {\bf{A $\Gamma$-loop}} at $\bar{a}$ is an arrow $a\rightarrow h.a$ and a {\bf{$\Gamma$-2 cycle}} at $\bar{a}$ is a pair of arrows $a\rightarrow j \rightarrow h.a$ such that $a,j \in Q_0$,  $j \notin \bar{a}$ and $h \in \Gamma$. \\
Let $Q$ be a locally-finite quiver with $B^\S$ {as its adjacency matrix and} a group of automorphisms {$\Gamma$} acting on it such that $Q$ does not admit a $\Gamma$-loop or a $\Gamma$-2 cycle at any of its orbits, the orbit-mutation in direction $\bar{k}$ is defined as follows \begin{center} \begin{equation}\label{eq1}\mu_{\bar{k}}(b^\S_{ij})=\begin{cases}
-b^\S_{ij}\hspace{3.1cm} \text{if}\hspace{0.1cm} i\in {\bar{k}} \hspace{0.1cm} \text{or}\hspace{0.1cm} j \in \bar{k}
\\
b^\S_{ij}+\underset{t \in \bar{k}}{\sum}\frac{\mid b^\S_{it}\mid b^\S_{tj}+b^\S_{it}\mid b^\S_{tj}\mid}{2}\hspace{0.3cm}\text{otherwise}\end{cases}
\end{equation}\end{center} Since $Q$ is locally-acyclic, {the summation in Equation (\ref{eq1})} is well-defined {and mutations in directions which belong to the same orbit commute since the quiver does not admit a $\Gamma-$loop}, hence we get the fact\begin{center}
\begin{equation}\label{ddd}\mu_{\bar{k}}(b^\S_{ij})=\underset{t \in {\bar{k}}\big\rvert_{\{i,j\}}}{\prod}\mu_t(b^\S_{ij})\end{equation}\end{center}
where ${\bar{k}}\big\rvert_{\{i,j\}}$ denotes the indices of $\bar{k}$ which are incident to $i$ or $j$ {and $\prod$ denotes the the composition of mutations in directions $t \in {\bar{k}}\big\rvert_{\{i,j\}}$}.\\

\begin{dfn}\label{coverunfoldingfolding}\begin{enumerate}\rm{
             \item \label{folding} Let $Q$ be a locally-finite quiver represented by $B^\S=(b_{ij}^\S)$ with no $ \Gamma$-loops or $\Gamma$-2 cycles and with the action of a group of automorphisms $\Gamma$ such that there are finitely many orbits $n<\infty$ under the action of this group. {The matrix} $B=(b_{\bar{i}\bar{j}})\in Mat_{n\times n}(\mathbb{Z})$ obtained by
$b_{\bar{i}\bar{j}}=\underset{k\in \bar{i}}{\sum}b^\S_{kj}$
is called the {\bf{folding}} of $Q$ and denoted by $B=B(Q)$.
             \item \label{cover}Conversely, let $B$ be a sing-skew-symmetric matrix such that there is a pair $(Q,\Gamma)$ where $Q$ is a (locally-finite) quiver             and $\Gamma$ is a group of automorphisms and $B=B(Q)$, then $(Q,\Gamma)$ is called a {\bf{covering}} of $B$.
             \item \label{unfolding12}If $(Q,\Gamma)$ is {a covering} of a sign-skew-symmetric matrix $B$ and $Q$ can perform arbitrary steps of orbit-mutation (the quiver obtained by any finite sequence of orbit-mutation does not have a $\Gamma$-loop or $\Gamma$-2 cycles), then $(Q,\Gamma)$ is called {an} {\bf{unfolding}} of $B$}.
           \end{enumerate}

           \end{dfn}
           \begin{re}
           {\rm{Through out this paper, sometimes we drop the group of automorphisms $\Gamma$ when pointing to an unfolding of a sign-skew-symmetric matrix and write $Q$ is an unfolding of $B$.}}
           \end{re}
{In \cite{mingli} the authors proved the following Lemma.}
\begin{lem}
\rm{Let $Q$ be a locally-finite quiver {and $\Gamma$} a group of automorphisms acting on it with finitely many number of orbits $\{\bar{i}_1, \bar{i}_2,....., \bar{i}_n\}$ such that $Q$ does not admit any $\Gamma$-loops or {$\Gamma$-2 cycles}, then the folding matrix $B$ of $Q$ is a sign-skew-symmetric matrix.}
\end{lem}
In what follows, we recall the construction that M. Huang and F. Li {set up} in \cite{mingli} to find a covering for acyclic sign-skew-symmetric matrices which can take arbitrary steps of orbit-mutation.
\begin{cons}\label{cons}
\rm{Let $B=(b_{ij}) \in Mat_{n \times n}(\mathbb{Z})$ be an acyclic sign-skew-symmetric matrix. An infinite quiver $Q(B)$ will be constructed inductively.
\begin{itemize}
  \item For each $i \in \{1,2,..,n\}$, we define a quiver $ Q^i$ as follows: $Q^i$ has\\ $\overset{n}{\underset{j=1}{\sum}}\mid b_{ji}\mid+1$ vertices with one vertex labeled by $i$ and other $\mid b_{ji} \mid$ {vertices} labeled by $j$ ($i\neq j$). If $b_{ji} > 0$ there is an arrow from each vertex labeled by $j$ to the unique vertex labeled by $i$. If $b_{ji} < 0$ there is an arrow from the unique vertex labeled by $i$ to each vertex labeled by $j$. No arrows between $i$ and $j$ if $b_{ij}=0$.
  \item Suppose we start the constructing process at $i=1$, we denote $Q^1=Q_{(1)}$. The unique vertex labeled by $1$ in $Q_{(1)}$ is called the old vertex, while the other vertices are called new vertices.
  \item For a new vertex in $Q_{(1)}$ labeled by $i_1$, $Q^{i_1}$ and $Q_{(1)}$ share a common arrow denoted by $\alpha_1$. We glue $Q_{({1})}$ and $Q^{i_1}$ along this common arrow. By iterating the gluing procedure for all $i_j \in I$ where $I$ is the set of the new vertices in $Q_{(1)}$, we get a new quiver $Q_{(2)}$ whose old vertices are the vertices of $Q_{(1)}$ and the other vertices are the new vertices. Clearly $Q_{(1)}$ is a subquiver of $Q_{(2)}$.
  \item Inductively, we obtain $Q_{({m+1})}$ from $Q_{({m})}$. Similarly, the old vertices are the vertices of $Q_{(m)}$ and the rest are new.
  \item Finally, we define the {(infinite)} quiver $Q(B)=\overset{\infty}{\underset{i=1}{\bigcup}}Q_{(i)}$, as $Q_{(m)}$ is always a subquiver of $Q_{({m+1})}$ for any $m$.
\end{itemize}
}
\end{cons}

\begin{re}\label{qi}
\rm{ Clearly we have the following facts:
\begin{enumerate}
\item The underlying quiver $Q(B)$ is a acyclic, since it is a tree clearly.

  \item The {full subquiver} of $Q(B)$ obtained by all the vertices incident to a vertex labeled by $i$ is $Q^i$.
  \item {Mostly, the quiver $Q(B)$ constructed as in Construction (\ref{cons}) is infinite but in some cases it might be finite. For example when $B$ is the adjacency skew-symmetric matrix of a finite tree $Q^\prime$, then $Q(B)=Q^\prime$ and thus $Q(B)$ is finite here.}

   \item Let $B^\S$ be the {(infinite)} skew-symmetric {matrix} corresponding to the {(infinite)} quiver $Q(B)$. The entries of $B^\S$ are either $-1,0$ or $1$.
  \item Let $\Gamma$ be a {subgroup
of the symmetric group $S_{Q(B)_0}$} that sends a vertex of $Q$ constructed as above to another vertex with the same label. By (2) in Remark \ref{qi}, the vertices which carry different labels are always connected to each other by the same way. That is if $b_{ij}^\S=a \in \{0,1,-1\}$, then $b_{g(i)g(j)}^\S=a$ for every $g \in \Gamma$ and hence $\Gamma$ is  the maximal subgroup of automorphisms which preserves the labels:\\{$\Gamma=\{h\in AutQ: \text{if}\hspace{0.1cm} h.a_s=h.a_t \hspace{0.1cm}\text{for}\hspace{0.1cm} a_s,a_t\in Q_0,\hspace{0.1cm} \text{then}\hspace{0.1cm} a_s,a_t\hspace{0.1cm} \text{have the same label}\}$}
\\By the action of $\Gamma$ all the vertices which have the same label lie in the same orbit.\\
M. Huang and F. Li in \cite{mingli}, proved the {following very important two }Theorems.\end{enumerate}
}
\end{re}
\begin{mytheo}\cite[Theorem 2.17]{mingli}
\rm{Any acyclic sign-skew-symmetric matrix $B$ of the size $n$ is always totally sign-skew-symmetric.}
\end{mytheo}
\begin{mytheo}\cite[Theorem 2.16]{mingli} \label{unfolding}
\rm{If $B$ is an {acyclic} sign-skew-symmetric matrix of the size $n \times n$, then $(Q(B),\Gamma)$ built from $B$ as in Construction \ref{cons} is an unfolding of $B$.}
\end{mytheo}

\begin{re}
{\rm{Through out the proof of Theorem \ref{unfolding} in \cite{mingli}, it was proved that the property of no $\Gamma-$ loops and no $\Gamma-$2 cycles is preserved under orbit-mutation for the (infinite) quiver $Q(B)$ constructed as in Construction \ref{cons} i.e, for any finite sequence of orbit-mutations the quiver $\mu_{\bar{k}_j}...\mu_{\bar{k}_1}(Q(B))$ does not admit any  $\Gamma-$ loops or $\Gamma-$2 cycles where $k_s\in Q_0(B)$ for every $1 \leq s\leq j$. This fact will be used later in this paper in places like the proof of Lemma \ref{orbitsource admissible}.}}
\end{re}
\section{The sign-coherence of $c-$vectors for an acyclic sign-skew-symmetric matrix}
\numberwithin{equation}{section}
{In this section, {we modify Construction \ref{cons} to find an unfolding of the $c-$vectores of} an extended sign-skew-symmetric matrix $\widetilde{B}=\begin{pmatrix}B\\I_n
 \end{pmatrix}\in Mat_{2n \times n}(\mathbb{Z})$.}\\
 Let $Q$ be a locally-finite quiver, the {{\bf{locally-finite framed quiver}}} $\widetilde{Q}$ {is the quiver} obtained from $Q$ by adding new vertices in a way that each vertex $a \in Q_0$ is connected to a new vertex $a^\prime$ by a single arrow $a\rightarrow a^\prime$ while $Q$ remains the same. {The elements of the set} $Q_0^\prime=\{a^\prime \mid a \in Q_0\}$ are called the {\bf{frozen vertices}}. This quiver is represented by the extended {infinite skew-symmetric} matrix $\widetilde{B^\S}=
\begin{pmatrix}B^\S\\I_\infty
\end{pmatrix}$. The bottom part of the matrix   $\widetilde{B^\S}$ is called the {\bf{$C$-matrix}} and $B^\S$ is called the {\bf{principal part}}.\\
{We extend the action of the group $\Gamma$ to the frozen vertices in the quiver $\widetilde{Q}$ such that {for every $g \in \Gamma$, $g(a^\prime)=g(b^\prime)$ if and only if $g(a)=g(b)$}, that is two frozen vertices lie in the same $\Gamma$-orbit if their mutable copies lie in the same orbit.} \\
{Let $\Gamma$ be a group of automorphisms acting on $B^\S$ such that $Q$ does not admit any $\Gamma-$loops or $\Gamma-$2 cycles, clearly $\Gamma$ is also a group of automorphisms of $I_\infty$. Hence $\Gamma$ is said to be a \textbf{group of automorphisms} of an extended matrix $\widetilde{B^\S}$ if it is a group of automorphisms of its principal part $B^\S$.}\\
We define the orbit-mutation on the $C$-matrix
 in direction $\bar{k}$ where $k\in Q_0$ as follows \begin{center} \begin{equation} \label{eq2}\mu_{\bar{k}}(c^\S_{ij})=\begin{cases}
-c^\S_{ij}\hspace{3.1cm}  \text{if}\hspace{0.1cm} j \in \bar{k}
\\
c^\S_{ij}+\underset{t \in \bar{k}}{\sum}\frac{\mid c^\S_{it}\mid b^\S_{tj}+c^\S_{it}\mid b^\S_{tj}\mid}{2}\hspace{0.3cm}\text{otherwise}\end{cases}
\end{equation}\end{center}
{Again since $Q$ does not admit a $\Gamma-$loop, the orbit-mutation of the $C-$matrix can be defined as
\begin{center}
\begin{equation}\label{rrr}
\mu_{\bar{k}}(c^\S_{ij})=\underset{t \in {\bar{k}}\big\rvert_{\{i,j\}}}{\prod}\mu_t(c^\S_{ij})
\end{equation}
\end{center}
 {Where $t \in {\bar{k}}\big\rvert_{\{i,j\}}$} denotes the indices of $\bar{k}$ which are incident to $i$ or $j$ {and $\prod$ denotes the the composition of mutations in directions $t \in {\bar{k}}\big\rvert_{\{i,j\}}$ }.\\
{By the definition of orbit-mutation for an extended infinite skew-symmetric matrix  $\widetilde{B^\S}=
\begin{pmatrix}B^\S\\I_\infty
\end{pmatrix}$ given in (\ref{eq1}) and (\ref{eq2}), it is easy to check that if $\Gamma$ is a group of automorphisms of $\widetilde{B^\S}$, it will be a group of automorphisms of any extended infinite skew-symmetric matrix obtained from $\widetilde{B^\S}$ by any finite sequence of orbit-mutations.}\\
{Since the extended adjacency matrix $\widetilde{B^\S}=\begin{pmatrix}B^\S\\I_\infty\end{pmatrix}$ of a locally-finite {framed} quiver $\widetilde{Q}$ is well-defined, we say that the \textbf{sign-coherence property} holds for a locally-finite {framed} quiver $\widetilde{Q}$ ({$\widetilde{Q}$ is sign-coherent}) if after performing any finite sequence of ordinary mutations {$\mu_{k_s}\mu_{k_{s-1}}...\mu_{k_1}$ on $\widetilde{B^\S}$ where $k_j \in Q_0$ for every $1\leq  j\leq s$}, the entries of any $c-$vector {in the matrix $\mu_{k_s}\mu_{k_{s-1}}...\mu_{k_1}(\widetilde{B^\S})$} are either all non-negative or all non-positive.}
{In other words a locally-finite framed quiver $\widetilde{Q}$ is sign-coherent if the arrows connecting any mutable vertex with the frozen vertices in the quiver $\mu_{k_s}\mu_{k_{s-1}}...\mu_{k_1}(\widetilde{Q})$} are either all emerging from this mutable vertex or all reaching at this mutable vertex, where $k_j \in Q_0$ for every $1\leq  j\leq s$ and $s< \infty$.}
{\begin{re}
\rm{The definitions of finite framed quivers and the sign-coherence property for finite framed quivers coincide with the definitions of these concepts for locally-finite quivers.}
\end{re}}

\begin{re}
{\rm{The quiver mutation for a framed quiver $\widetilde{Q}$ (finite or locally-finite) can be taken only in direction of a mutable vertex and is obtained as for the quiver $Q$ with one modification which is to remove all arrows that can be created during the mutation process between any two frozen vertices.}}
\end{re}

The sign-coherence property was proved for {finite skew-symmetric matrices} (finite quivers) in \cite{Derksen} and for {finite skew-symmetrizable matrices} in \cite{GR}. Here we prove that {the sign-coherence property holds} for locally-finite {framed} quivers.
 \begin{lem}\label{infinitesigncoherence}
 \rm{The sign-coherence property {of $c-$vectors} holds for locally-finite {framed} quivers.}
 \begin{proof}
 {Suppose that $\widetilde{Q}$ is a locally-finite framed quiver, and $\mu_{k_m}......\mu_{k_1}$ is a composition of mutations such that $  k_j \in Q_0$ for every $1 \leq j \leq m$. We consider the whole quiver $\widetilde{Q}$ as a quiver resulting from the gluing of two {full} subquivers of $\widetilde{Q}$ as follows
 \begin{itemize}
                                            \item The first full subquiver of $\widetilde{Q}$ is $\widetilde{Q}\big\rvert_{k_m..k_1}$ which is obtained by the vertices $\{k_m,k_{m-1},..,k_1\}$ and their frozen copies and the mutable vertices in $N$ with their frozen copies such that $N$ is the set of mutable vertices that are in the neighbourhood of $\{k_m,k_{m-1},..,k_1\}$. Clearly $\widetilde{Q}\big\rvert_{k_m..k_1}$ is a finite framed quiver.
                                            \item The other full subquiver of $\widetilde{Q}$ is $\widetilde{Q}\big\rvert_{S}$ which is obtained by the vertices of $S$ and their frozen copies where $S$ is the set of mutable vertices in $\widetilde{Q}$ which are not contained in $\widetilde{Q}\big\rvert_{k_m..k_1}$. Clearly $\widetilde{Q}\big\rvert_{S}$ is a locally-finite framed quiver satisfying the sign-coherence property by construction.
  \end{itemize}
The gluing procedure occurs between $\widetilde{Q}\big\rvert_{k_m..k_1}$ and $\widetilde{Q}\big\rvert_{S}$ along the set of arrows $A$ in $\widetilde{Q}$ connecting $N$ and $S$.\\
 Since mutation at some mutable vertex $r$ reverses the direction of all arrows incident to the vertex $r$ and may affect the arrows between the vertices in the neighbourhood of $r$, the quiver $\widetilde{Q}\big\rvert_{S}$ and the set of arrows $A$ remain unchanged during the composition of mutations $\mu_{k_m}......\mu_{k_1}$.\\
Hence $\mu_{k_m}.....\mu_{k_1}(\widetilde{Q}\big\rvert_{S})=\widetilde{Q}\big\rvert_{S}$ and the quiver $\mu_{k_m}.....\mu_{k_1}(\widetilde{Q})$ is the gluing of two quivers along the set of arrows $A$ described as above and which connects mutable vertices. These two quivers are :
\begin{itemize}
  \item The first one is $\mu_{k_m}.....\mu_{k_1}(\widetilde{Q}\big\rvert_{k_m..k_1})$ which is sign-coherent for $\widetilde{Q}\big\rvert_{k_m..k_1}$ is a finite framed-quiver \cite{Derksen}, \cite{GR}.
  \item The other one is $\mu_{k_m}.....\mu_{k_1}(\widetilde{Q}\big\rvert_{S})=\widetilde{Q}\big\rvert_{S}$ which is sign-coherent by construction.
\end{itemize}
Thus $\widetilde{Q}$ is sign-coherent.}

\end{proof}
 \end{lem}
 {Let $\widetilde{A}=\begin{pmatrix}
 A\\C\end{pmatrix}$ be a matrix with $C$ as the $C-$matrix and its columns are the $c-$vectors. Suppose that this matrix is equipped with a group of automorphisms $\Gamma$ acting on its principal part $A$ such that the ordinary mutation and orbit-mutation are well-defined and $ \widetilde{A}$ can do arbitrary steps of ordinary mutation and orbit-mutation.
 If all the entries of any $c-$vector are either all non-positive or all non-negative after performing any finite sequence of orbit-mutation i.e, any $c-$vector in the matrix $\mu_{\bar{i}_m}....\mu_{\bar{i}_0}(\widetilde{A})$ has entries which are all non-negative or all non-positive  where $i_j$ is used to index the matrix $A$ for $0 \leq j \leq m$ and $0 \leq m < \infty$, then we say that the {\bf{orbit-sign coherence property}} holds for this matrix.} By convention $\mu_{\bar{i}_0}(\widetilde{A})=\widetilde{A}$.({{$\mu_{\bar{i}_0}$} means no orbit-mutation has been applied yet and any $c-$vector in $\widetilde{A}$ has its entries either all non-positive or all non-negative}.)

 \begin{cor}\label{infiniteorbitsign}
 {\rm{ The orbit-sign coherence property holds for a locally-finite framed quiver {$\widetilde{Q}$} equipped with a group of {automorphisms} $\Gamma$ and can do arbitrary steps of orbit-mutation.}}
 \begin{proof}
 {For a locally-finite {framed} quiver}, any finite sequence of orbit mutation can be regarded as a longer but still finite sequence of ordinary mutation (see Equations (\ref{ddd}) and (\ref{rrr})) and by Lemma \ref{infinitesigncoherence} the entries of any $c-$vector of the quiver obtained by any finite sequence of ordinary mutation performed on a locally-finite {framed} quiver $\widetilde{Q}$ are either all non-negative or all non-positive. Hence the result follows.
 \end{proof}
\end{cor}
{We denote by $sgn(i)$ the sign of the  column indexed by $i$ in the $C-$matrix and $sgn(i)=+$ when the entries of the column $i$ are non-negative while $sgn(i)=-$ when the entries of the column $i$ are non-positive.}
\begin{lem}
\label{orbitlemmalocallyfinite}
 \rm{Let $\widetilde{Q}$ be a locally-finite {framed} quiver equipped with a group of automorphisms $\Gamma$ such that $\widetilde{Q}$ does not admit a $\Gamma-$loop or $\Gamma-$2 cycle with $\widetilde{B^\S}=
\begin{pmatrix}B^\S\\I_\infty
\end{pmatrix}$ as its adjacency matrix, and let ${\widetilde{B^\S}}^{\bar{k}}=\mu_{\bar{k}}(\widetilde{B^\S})=\begin{pmatrix}
{B^\S}^{\bar{k}}\\{C^\S}^{\bar{k}}\end{pmatrix}$ be the matrix obtained by orbit-mutation in direction $\bar{k}$ such that $ k \in Q_0$. If $i_1,i_2$ {fall in the} same orbit, the columns indexed by $i_1, i_2  $ in the $C$-matrix ${C^\S}^{\bar{k}}$ have the same sign.}
\begin{proof}
{Since $i_1$ and $i_2$ fall in the same orbit, there exists an automorphism $g\in \Gamma$ such that $i_2=g(i_1)$. Clearly $\Gamma$ is a group of automorphisms of ${\widetilde{B^\S}}^{\bar{k}}=\mu_{\bar{k}}(\widetilde{B^\S})=\begin{pmatrix}
{B^\S}^{\bar{k}}\\{C^\S}^{\bar{k}}\end{pmatrix}$. Thus ${{c^\S}^{\bar{k}}}_{li_1}={{{c^{\S}}^{\bar{k}}_{g(l)g(i_1)}}}={{{c^{\S}}^{\bar{k}}_{g(l)i_2}}}$ for any index $l$, thus $sgn(i_1)=sgn(i_2)$.}
\end{proof}
\end{lem}
{When $sgn(s)=+$ for every $s\in \bar{i}$, then $sgn(\bar{i})=+$ and $\bar{i}$ is said to be a {\bf{green orbit}}. Respectively, when $sgn(s)=-$ for every $s\in \bar{i}$, then $sgn(\bar{i})=-$ and $\bar{i}$ is said to be a {\bf{red orbit}}.}

Let ${\widetilde{B^\S}}=\begin{pmatrix}
{B^\S}\\{C^\S}\end{pmatrix}$ be the adjacency matrix of a locally-finite {framed} quiver $\widetilde{Q}$, in \cite{mingli} the authors defined the folding matrix $B=(b_{\bar{i}\bar{j}})$ of a locally-finite quiver $Q$ endowed with a group of automorphisms $\Gamma$ \begin{center}
$b_{\bar{i}\bar{j}}=\underset{k\in \bar{i}}{\sum}b^\S_{kj}$\end{center}
{Analogously we define the folding of the $C$-matrix} \begin{center}
$c_{\bar{i}\bar{j}}=\underset{k\in \bar{i}}{\sum}c^\S_{kj}$\end{center}
Clearly, when we have a finite number $n$ of orbits, the folding of the adjacency matrix $\widetilde{B^\S}$ of a locally-finite {framed} quiver $\widetilde{Q}$, is the extended matrix $\widetilde{B}=\begin{pmatrix}
B\\I_n
\end{pmatrix} \in Mat_{2n\times n}(\mathbb{Z})$.

            Now we construct an unfolding for a given extended sign-skew-symmetric matrix.
\begin{cons}\label{cons2}
\rm{Let $\widetilde{B}=\begin{pmatrix}B\\I_n
 \end{pmatrix}\in Mat_{2n \times n}(\mathbb{Z})$ be an extended acyclic sign-skew-symmetric matrix. A {(locally-finite)} framed quiver {$\widetilde{Q}(\widetilde{B})$} will be constructed inductively.
 \begin{itemize}
  \item For each mutable vertex $i \in \{1,2,..,n\}$, we define a quiver $ \widetilde{Q^i}$ as follows: $\widetilde{Q^i}$ has $\overset{n}{\underset{j=1}{\sum}}\mid b_{ji}\mid+2$ vertices with one vertex labeled by $i$ and other $\mid b_{ji} \mid$ {vertices} labeled by $j$ ($i\neq j$). If {$b_{ji} < 0$} there is an arrow from each vertex labeled by $j$ to the unique vertex labeled by $i$. If {$b_{ji} >0$} there is an arrow from the unique vertex labeled by $i$ to each vertex labeled by $j$. No arrows between $i$ and $j$ if $b_{ij}=0$. And finally with one vertex labeled by $i^\prime$ which is the frozen copy of $i$ such that there is one arrow $i\rightarrow i^\prime$.
      \item {We start by considering  $ \widetilde{Q^1}$ as the initial subquiver and we denote $\widetilde{Q_{(1)}}=\widetilde{Q^1}$. During the constructing process which has $\widetilde{Q_{(1)}}$ as its initial subquiver, the mutable vertices are either old or new while the frozen vertices are not considered old or new. For $\widetilde{Q_{(1)}}$ the vertex $1$ is an old vertex and the other mutable vertices are new. For every new vertex $i$, $\widetilde{Q^i}$ and $\widetilde{Q_{(1)}}$ share a common arrow $\alpha_i$, we glue $\widetilde{Q^i}$ and $\widetilde{Q_{(1)}}$ along this common arrow to get a new subquiver $\widetilde{Q_{(2)}}$. The old vertices of $\widetilde{Q_{(2)}}$ are the mutable vertices of $\widetilde{Q_{(1)}}$ and the other mutable vertices are new.}
          \item We continue inductively as in Construction \ref{cons} and build $\widetilde{Q_{(m+1)}}$ from $\widetilde{Q_{(m)}}$.
           \item We define {$\widetilde{Q}(\widetilde{B})=\overset{\infty}{\underset{i=1}{\bigcup}}\widetilde{Q_{(i)}}$}.
          \end{itemize}}
\end{cons}
\begin{re}\label{unfolding2}
\rm{{
\begin{enumerate}
  \item The matrix associated with the (locally-finite) quiver $\widetilde{Q}(\widetilde{B})$ obtained from Construction \ref{cons2} is the (infinite) and well-defined matrix $\widetilde{B^\S}=\begin{pmatrix}B^\S\\I_{\infty}\end{pmatrix} $  where the upper part of this matrix is the principal part such that {$b^\S_{ij}<0$} if there are $\mid b^\S_{ij}\mid$ many arrows from the mutable vertex $i$ to the mutable vertex $j$ whereas {$b^\S_{ij}>0$} if there are $\mid b^\S_{ij}\mid$ many arrows from the mutable vertex $j$ to the mutable vertex $i$ and $b^\S_{ij}=0$ if there are no arrows between the mutable vertices $i$ and $j$. {The lower part of this matrix is the $C-$ matrix such that $c^\S_{ij}>0$ if there are $\mid c^\S_{ij}\mid$ many arrows from the mutable vertex $j$ to the frozen vertex $i^\prime$ whereas $c^\S_{ij}<0$ if there are $\mid c^\S_{ij}\mid$ many arrows from the frozen vertex $i^\prime$ to the mutable vertex $j$ and $c^\S_{ij}=0$ if there are no arrows between the mutable vertex $j$ and the frozen vertex $i^\prime$.}
  \item Let $\Gamma$ be the maximum subgroup that preserves labels of the symmetric matrix $S_{{\widetilde{Q}(\widetilde{B})}_0}$ acting on the set of mutable and frozen vertices of $\widetilde{Q}(\widetilde{B})$. By Construction \ref{cons2}, $\Gamma$ is a group of automorphisms and the orbits obtained by its action are $\{\bar{1},...\bar{n}, \bar{1}^\prime,...,\bar{n}^\prime\}$ such that the orbit $\bar{i}$ contains all the mutable vertices labeled by $i$ and the orbit $\bar{i}^\prime$ contains all the frozen vertices labeled by $i^\prime$ for every $1\leq i \leq n$.
  \item Clearly, the folding of the adjacency matrix $\widetilde{B^\S}$ of the quiver $\widetilde{Q}(\widetilde{B})$ is $\widetilde{B}$. The full subquiver $Q$ of $\widetilde{Q}$ obtained by the mutable vertices is exactly as constructed in Construction \ref{cons} {taking into consideration the way we follow in this paper to associate a quiver with a matrix}  thus it is an unfolding of $B$ by Theorem \ref{unfolding} that is $Q$ can take arbitrary steps of orbit mutations and since orbit-mutation is only taken in direction of an orbit whose elements are labels of mutable vertices, we conclude that $\widetilde{Q}(\widetilde{B})$ is an unfolding of $\widetilde{B}$.
\end{enumerate}
}}
\end{re}

           \begin{re}
{\rm{Let $\widetilde{B}$ be the folding of $\widetilde{B^\S}$ associated {with} a quiver $\widetilde{Q}$, to avoid ambiguity, when we mutate {$\widetilde{B}$} of in direction $\bar{i}$, the mutation will be denoted as $\mu_{i^f}$ {since this mutation is an ordinary mutation here and not orbit-mutation}.}}
\end{re}
            By convection $\mu^{\bar{\sigma}_0}(\widetilde{Q})=\widetilde{Q}$ and  $\mu^{\sigma_0^f}(\widetilde{B})=\widetilde{B}$.
           \begin{lem}\label{invariantunfolding}
          \rm{ Let $\widetilde{B}=\begin{pmatrix}
          B\\I_n
          \end{pmatrix}$ be an extended, finite and acyclic sign-skew-symmetric matrix and let $\widetilde{Q}(\widetilde{B})$ be the {(locally-finite) {framed} quiver associated with the (infinite) and well-defined skew-symmetric matrix $\widetilde{B^{\S}}=\begin{pmatrix}B^{\S}\\I_{\infty}
          \end{pmatrix}$ obtained from Construction \ref{cons2} as an unfolding of $\widetilde{B}$.
           We denote by $\mu^{\bar{\sigma}_m}(\widetilde{B^{\S}})$ the composition of orbit mutations {$\mu_{\bar{k}_m}....\mu_{\bar{k}_0}(\widetilde{B^{\S}})$} such that $k_j \in Q_0$ for {$0\leq j \leq m$} and we denote by $\mu^{\sigma_m^f}(\widetilde{B})$ the composition of ordinary mutation $\mu_{k_m^f}....\mu_{k_1^f}(\widetilde{B})$. Then $\mu^{\bar{\sigma}_m}(\widetilde{B^\S})$ is {an} unfolding of $\mu^{\sigma_m^f}(\widetilde{B})$.}}
          \begin{proof}
          Since $\widetilde{Q}(\widetilde{B})$ can take arbitrary steps of orbit-mutation, $\mu^{\bar{\sigma}_m}(\widetilde{Q}(\widetilde{B}))$ can also take arbitrary steps of orbit-mutation {so} we need only to prove that $\mu^{\bar{\sigma}_m}(\widetilde{Q}(\widetilde{B}))$ is {a covering} of $\mu^{\sigma_m^f}(\widetilde{B})$. We prove it by induction for the mutable part $Q(B)$ {represented by the matrix $B^\S$} first, then for the frozen part represented by the $C$-matrix.  Trivially, $\mu^{\bar{\sigma}_0}(Q(B))$ is {a covering} of $\mu^{\sigma_0^f}(B)$. Suppose the result holds for every $v < m$.

           $\mu^{\bar{\sigma}_m}(b^\S_{ij})=\begin{cases}
              -\mu^{\bar{\sigma}_{m-1}}(b^\S_{ij})\hspace{7.1cm} \text{if}\hspace{0.1cm} i\in \bar{k}_m \hspace{0.1cm} \text{or}\hspace{0.1cm} j \in \bar{k}_m
\\
               \mu^{\bar{\sigma}_{m-1}}(b^\S_{ij})+\underset{t \in \bar{k}_m}{\sum}\frac{\mid \mu^{\bar{\sigma}_{m-1}}(b^\S_{it})\mid \mu^{\bar{\sigma}_{m-1}}( b^\S_{tj})+\mu^{\bar{\sigma}_{m-1}}(b^\S_{it})\mid \mu^{\bar{\sigma}_{m-1}}( b^\S_{tj})\mid}{2}\hspace{0.3cm}\text{otherwise}\end{cases}$\\ When $i \in \bar{k}_m$ or $j \in \bar{k}_m$, $\underset{s \in \bar{i}}{\sum}\mu^{\bar{\sigma}_m}(b_{sj}^\S)=-\underset{s \in \bar{i}}{\sum}\mu^{\bar{\sigma}_{m-1}}(b_{sj}^\S)=-\mu^{\sigma_{m-1}^f}(b_{\bar{i}\bar{j}})=\mu^{{\sigma_{m}^f}}(b_{\bar{i}\bar{j}})$.\\\\ When $i \notin \bar{k}_m$ and $j \notin \bar{k}_m$,\\\\ $\underset{s \in \bar{i}}{\sum}\mu^{\bar{\sigma}_m}(b_{sj}^\S)=\underset{s \in \bar{i}}{\sum}\bigl(\mu^{\bar{\sigma}_{m-1}}(b^\S_{sj})+\underset{t \in \bar{k}_m}{\sum}\frac{\mid \mu^{\bar{\sigma}_{m-1}}(b^\S_{st})\mid \mu^{\bar{\sigma}_{m-1}}( b^\S_{tj})+\mu^{\bar{\sigma}_{m-1}}(b^\S_{st})\mid \mu^{\bar{\sigma}_{m-1}}( b^\S_{tj})\mid}{2}\bigr)=\\\\\underset{s \in \bar{i}}{\sum}\bigl(\mu^{\bar{\sigma}_{m-1}}(b^\S_{sj})\bigr)+\underset{s \in \bar{i}}{\sum}\bigl(\underset{t \in \bar{k}_m}{\sum}\frac{\mid \mu^{\bar{\sigma}_{m-1}}(b^\S_{st})\mid \mu^{\bar{\sigma}_{m-1}}( b^\S_{tj})+\mu^{\bar{\sigma}_{m-1}}(b^\S_{st})\mid \mu^{\bar{\sigma}_{m-1}}( b^\S_{tj})\mid}{2}\bigr)=\\\\\mu^{{\sigma}^f_{m-1}}(b_{\bar{i}\bar{j}})+\underset{t \in \bar{k}_m}{\sum}\bigl(\underset{s \in \bar{i}}{\sum}\frac{\mid \mu^{\bar{\sigma}_{m-1}}(b^\S_{st})\mid \mu^{\bar{\sigma}_{m-1}}( b^\S_{tj})+\mu^{\bar{\sigma}_{m-1}}(b^\S_{st})\mid \mu^{\bar{\sigma}_{m-1}}( b^\S_{tj})\mid}{2}\bigr)$\\\\
                Since $\mu^{\bar{\sigma}_{m-1}}(Q(B))$ does not have a $\Gamma$-2 cycles, when $s \in \bar{i}$ all the entries $\mu^{\bar{\sigma}_{m-1}}(b_{st}^\S)$ have the same sign and when $t \in \bar{k}_m$, the entries $\mu^{\bar{\sigma}_{m-1}}(b_{tj}^\S)$ have the same sign, hence\\\\
                $\underset{s \in \bar{i}}{\sum}\mu^{\bar{k}_m}(b_{sj}^\S)=\mu^{{\sigma}^f_{m-1}}(b_{\bar{i}\bar{j}})+\underset{t \in \bar{k}_m}{\sum}\bigl(\frac{\mid \underset{s \in \bar{i}}{\sum}\mu^{\bar{\sigma}_{m-1}}(b^\S_{st})\mid \mu^{\bar{\sigma}_{m-1}}( b^\S_{tj})+\underset{s \in \bar{i}}{\sum}\mu^{\bar{\sigma}_{m-1}}(b^\S_{st})\mid  \mu^{\bar{\sigma}_{m-1}}( b^\S_{tj})   \mid}{2}$\\\\
                $\mu^{{\sigma}^f_{m-1}}(b_{\bar{i}\bar{j}})+\frac{\mid \mu^{{\sigma}^f_{m-1}}(b_{\bar{i}\bar{t}})\mid \underset{t \in \bar{k}_m}{\sum}( \mu^{\bar{\sigma}_{m-1}}(b^\S_{tj}))+\mu^{{\sigma}^f_{m-1}}(b_{\bar{i}\bar{t}})\mid  \underset{t \in \bar{k}_m}{\sum}(\mu^{\bar{\sigma}_{m-1}}( b^\S_{tj}))   \mid}{2}$\\\\
                $=\mu^{{\sigma}^f_{m-1}}(b_{\bar{i}\bar{j}})+\frac{\mid \mu^{{\sigma}^f_{m-1}}(b_{\bar{i}\bar{t}})\mid   \mu^{{\sigma}^f_{m-1}}( b_{\bar{t}\bar{j}})+\mu^{{\sigma}^f_{m-1}}(b_{\bar{s}\bar{t}})\mid  \mu^{{\sigma}^f_{m-1}}( b_{\bar{t}\bar{j}})   \mid}{2}=\mu^{{\sigma}^f_{m}}( b_{\bar{i}\bar{j}})$.
                And now we will prove that {a covering} of the $C-$matrix is invariant under a composition of orbit mutation. Trivially,$\mu^{\bar{\sigma}_0}(I_\infty)$ is a covering of $\mu^{\sigma_0^f}(I_n)$.\\
                $\mu^{\bar{\sigma}_m}(c^\S_{ij})=\begin{cases}
              -\mu^{\bar{\sigma}_{m-1}}(c^\S_{ij})\hspace{7.1cm} \text{if}\hspace{0.1cm} j \in \bar{k}_m
\\
               \mu^{\bar{\sigma}_{m-1}}(c^\S_{ij})+\underset{t \in \bar{k}_m}{\sum}\frac{\mid \mu^{\bar{\sigma}_{m-1}}(c^\S_{it})\mid \mu^{\bar{\sigma}_{m-1}}( b^\S_{tj})+\mu^{\bar{\sigma}_{m-1}}(c^\S_{it})\mid \mu^{\bar{\sigma}_{m-1}}( b^\S_{tj})\mid}{2}\hspace{0.3cm}\text{otherwise}\end{cases}$\\ When  $j \in \bar{k}_m$, $\underset{s \in \bar{i}}{\sum}\mu^{\bar{\sigma}_m}(c_{sj}^\S)=-\underset{s \in \bar{i}}{\sum}\mu^{\bar{\sigma}_{m-1}}(c_{sj}^\S)=-\mu^{\sigma_{m-1}^f}(c_{\bar{i}\bar{j}})=\mu^{{\sigma_{m}^f}}(c_{\bar{i}\bar{j}})$.\\\\ When $j \notin \bar{k}_m$,\\\\ $\underset{s \in \bar{i}}{\sum}\mu^{\bar{\sigma}_m}(c_{sj}^\S)=\underset{s \in \bar{i}}{\sum}\bigl(\mu^{\bar{\sigma}_{m-1}}(c^\S_{sj})+\underset{t \in \bar{k}_m}{\sum}\frac{\mid \mu^{\bar{\sigma}_{m-1}}(c^\S_{st})\mid \mu^{\bar{\sigma}_{m-1}}( b^\S_{tj})+\mu^{\bar{\sigma}_{m-1}}(c^\S_{st})\mid \mu^{\bar{\sigma}_{m-1}}( b^\S_{tj})\mid}{2}\bigr)=\\\\\underset{s \in \bar{i}}{\sum}\bigl(\mu^{\bar{\sigma}_{m-1}}(c^\S_{sj})\bigr)+\underset{s \in \bar{i}}{\sum}\bigl(\underset{t \in \bar{k}_m}{\sum}\frac{\mid \mu^{\bar{\sigma}_{m-1}}(c^\S_{st})\mid \mu^{\bar{\sigma}_{m-1}}( b^\S_{tj})+\mu^{\bar{\sigma}_{m-1}}(c^\S_{st})\mid \mu^{\bar{\sigma}_{m-1}}( b^\S_{tj})\mid}{2}\bigr)=\\\\\mu^{{\sigma}^f_{m-1}}(c_{\bar{i}\bar{j}})+\underset{t \in \bar{k}_m}{\sum}\bigl(\underset{s \in \bar{i}}{\sum}\frac{\mid \mu^{\bar{\sigma}_{m-1}}(c^\S_{st})\mid \mu^{\bar{\sigma}_{m-1}}( b^\S_{tj})+\mu^{\bar{\sigma}_{m-1}}(c^\S_{st})\mid \mu^{\bar{\sigma}_{m-1}}( b^\S_{tj})\mid}{2}\bigr)$\\\\
               By Corollary \ref{infiniteorbitsign}, when $s \in \bar{i}$ all the entries $\mu^{\bar{\sigma}_{m-1}}(c_{st}^\S)$ have the same sign and hence\\\\
                $\underset{s \in \bar{i}}{\sum}\mu^{\bar{\sigma}_m}(c_{sj}^\S)=\mu^{{\sigma}^f_{m-1}}(c_{\bar{i}\bar{j}})+\underset{t \in \bar{k}_m}{\sum}\bigl(\frac{\mid \underset{s \in \bar{i}}{\sum}\mu^{\bar{\sigma}_{m-1}}(c^\S_{st})\mid \mu^{\bar{\sigma}_{m-1}}( b^\S_{tj})+\underset{s \in \bar{i}}{\sum}\mu^{\bar{\sigma}_{m-1}}(c^\S_{st})\mid  \mu^{\bar{\sigma}_{m-1}}( b^\S_{tj})   \mid}{2}$\\\\
                $\mu^{{\sigma}^f_{m-1}}(c_{\bar{i}\bar{j}})+\frac{\mid \mu^{{\sigma}^f_{m-1}}(c_{\bar{i}\bar{t}})\mid \underset{t \in \bar{k}_m}{\sum}( \mu^{\bar{\sigma}_{m-1}}(b^\S_{tj}))+\mu^{{\sigma}^f_{m-1}}(c_{\bar{i}\bar{t}})\mid  \underset{t \in \bar{k}_m}{\sum}(\mu^{\bar{\sigma}_{m-1}}( b^\S_{tj}))   \mid}{2}$\\\\
                $=\mu^{{\sigma}^f_{m-1}}(c_{\bar{i}\bar{j}})+\frac{\mid \mu^{{\sigma}^f_{m-1}}(c_{\bar{i}\bar{t}})\mid   \mu^{{\sigma}^f_{m-1}}( b_{\bar{t}\bar{j}})+\mu^{{\sigma}^f_{m-1}}(c_{\bar{s}\bar{t}})\mid  \mu^{{\sigma}^f_{m-1}}( b_{\bar{t}\bar{j}})   \mid}{2}=\mu^{{\sigma}^f_{m}}( c_{\bar{i}\bar{j}})$.
               \end{proof}

           \end{lem}

\begin{eg}
\rm{The construction of an unfolding of the extended acyclic sign-skew symmetric matrix {\begin{center} $\widetilde{B}=\begin{pmatrix}
0 & -1 & 0 & -1\\
3 & 0 & -1 & 0\\
0 & 5 & 0 & -2\\
1 & 0 & 3 & 0\\
&...................&\\
1& 0 & 0 & 0\\
0 & 1 & 0 & 0\\
0 & 0 & 1 & 0\\
0 & 0 & 0 & 1
\end{pmatrix}$\end{center}}
built according to Construction \ref{cons2} is shown in Figrue \ref{coffee} and Figure \ref{boat}}.
\end{eg}

\begin{figure}[h!]
  \centering
  \begin{subfigure}[b]{0.2\linewidth}
    \includegraphics[width=\linewidth]{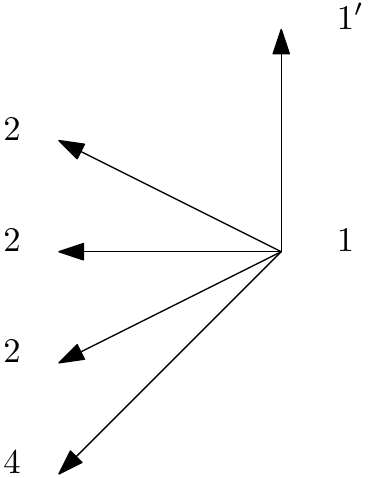}
    \caption{$\widetilde{Q^1}$}
  \end{subfigure}
  \begin{subfigure}[b]{0.2\linewidth}
    \includegraphics[width=\linewidth]{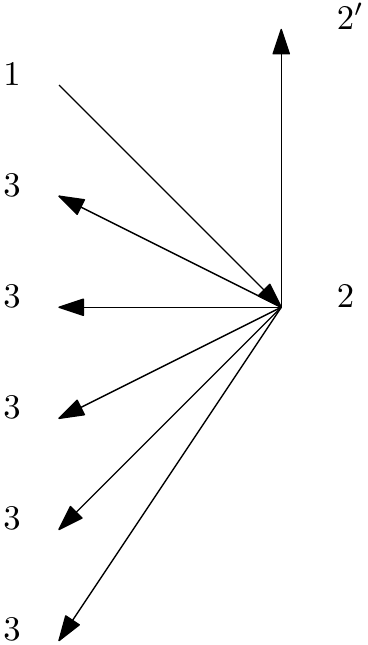}
    \caption{$\widetilde{Q^2}$}
  \end{subfigure}
  \begin{subfigure}[b]{0.2\linewidth}
    \includegraphics[width=\linewidth]{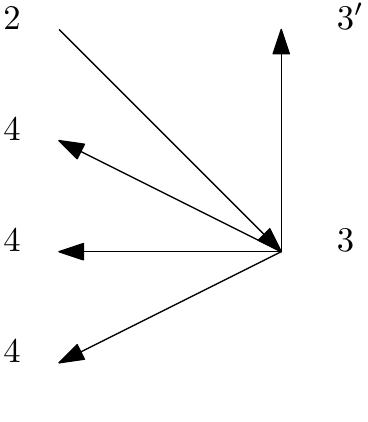}
    \caption{$\widetilde{Q^3}$}
  \end{subfigure}
  \begin{subfigure}[b]{0.2\linewidth}
    \includegraphics[width=\linewidth]{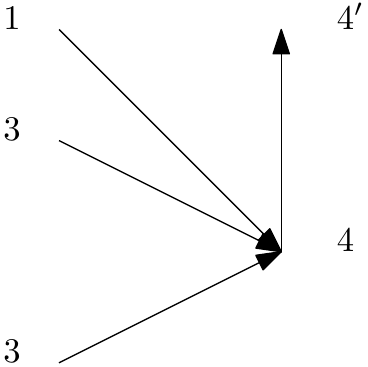}
    \caption{$\widetilde{Q^4}$}
  \end{subfigure}
  \caption{}
  \label{coffee}
\end{figure}

\begin{figure}[h!]
  \includegraphics[width=\linewidth]{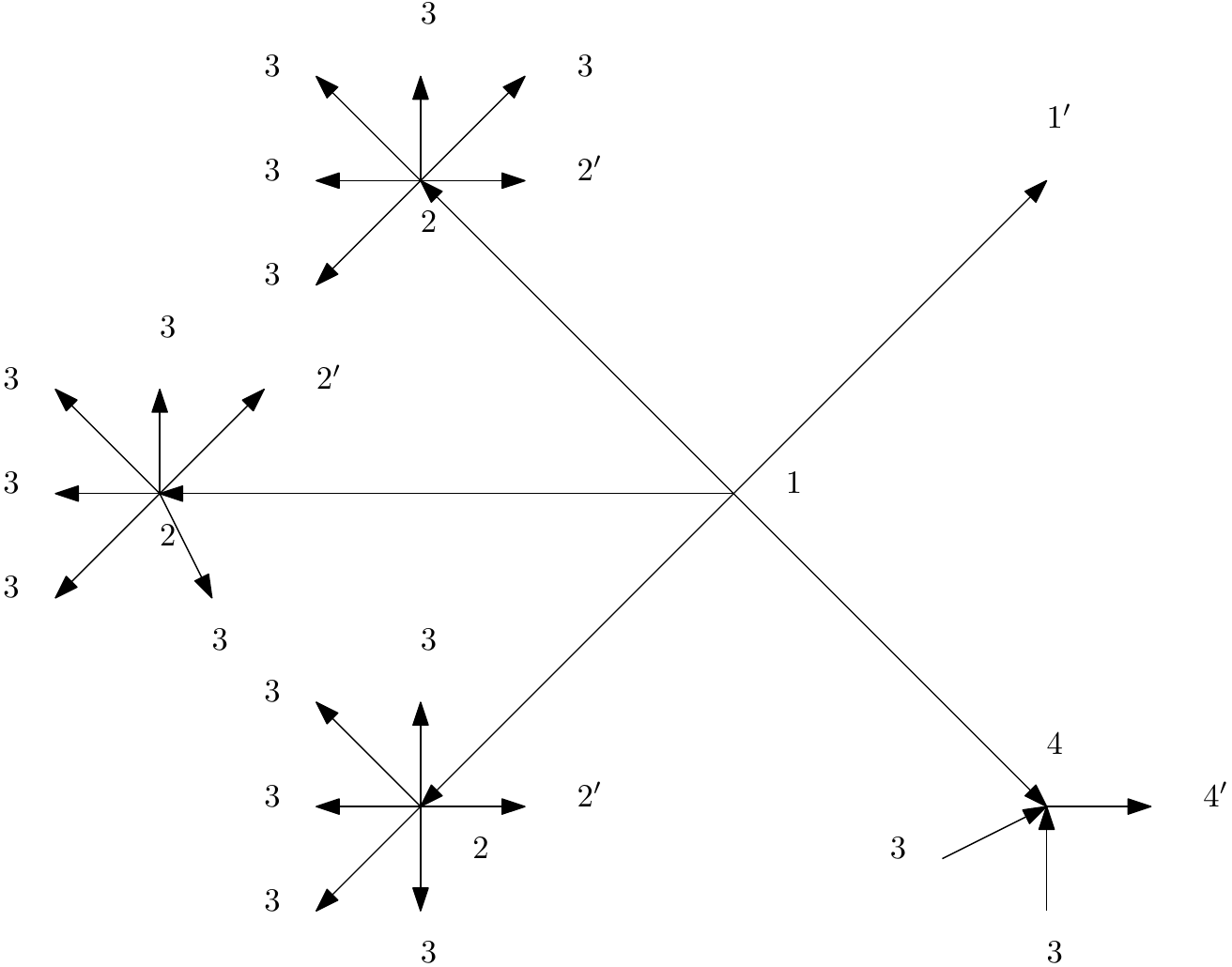}
  \caption{$\widetilde{Q_{(1)}}$}
  \label{boat}
\end{figure}
\begin{mytheo}
\rm{The sign-coherence property of {$c-$vectors} holds for acyclic sign-skew-symmetric matrices.}
\begin{proof}
{Let $\widetilde{B}=\begin{pmatrix}B\\I_n\end{pmatrix} \in Mat_{2n \times n}{\mathbb{Z}}$ be an extended sign-skew-symmetric matrix and let {$\widetilde{Q}(\widetilde{B})$} be the locally-finite {framed} quiver obtained from Construction \ref{cons2} as an unfolding of $\widetilde{B}$ with the adjacency matrix $ \widetilde{B^\S}=\begin{pmatrix}B^\S\\I_\infty\end{pmatrix}$. We denote by $\mu^{\bar{\sigma}_k}(\widetilde{B^\S})=\begin{pmatrix}
{B^\S}^{\bar{\sigma}_k}\\
{C^\S}^{{\bar{\sigma}_k}}
\end{pmatrix}$  the matrix obtained from $\begin{pmatrix}
{B^\S}\\
I_\infty
\end{pmatrix}$ after taking a finite sequence $(\mu_{\bar{i}_0},...,\mu_{\bar{i}_k})$ of orbit-mutation. We denote by $\mu^{{\sigma}_k^f}(\widetilde{B})=\begin{pmatrix}
{B}^{{\sigma}_k^f}\\
{C}^{{{\sigma}_k^f}}
\end{pmatrix}$ the matrix obtained from $\begin{pmatrix}
{B}\\
I_n
\end{pmatrix}$ after taking a finite sequence $(\mu_{{i_0}}^f,...,\mu_{{i_k}}^f)$ of {ordinary} mutation such that $i_j^f$ refers to the order of the row or column indexed by $\bar{i}_j$ in the folding matrix where $ i_j \in Q_0$ for $0 \leq j \leq k$. By convention $ \mu^{\bar{\sigma}_{_0}}(\widetilde{B^\S})=\widetilde{B^\S}$ and $\mu^{{\sigma}_{0}^f}(\widetilde{B})=\widetilde{B}$. By Lemma \ref{invariantunfolding}  $\mu^{\bar{\sigma}_{_j}}(\widetilde{B^\S})$ is {an} unfolding of $\mu^{\sigma_{j}^f}(\widetilde{B})$ for $0 \leq j\leq k$. {{The orbit-sign coherence property holds for the locally-finite framed quiver $\widetilde{Q}(\widetilde{B})$ by Lemma \ref{infiniteorbitsign}}}, {thus} the entries of any $c-$vector in the matrix $\mu^{\bar{\sigma}_k}(\widetilde{B^\S})$ are either all non-negative or all non-positive. By the definition of {an} unfolding, we find that the entries of any $c-$vector in ${B}^{{{\sigma}_k^f}}$ are either all non-negative or non-positive for for any $0 \leq k < \infty$ and thus the sign-coherence property holds for $\widetilde{B}$.}

\end{proof}

\end{mytheo}

\section{Maximal green sequences for an acyclic sign-skew-symmetric matrix}
After {proving} that the sign-coherence property holds for an acyclic sign-skew-symmetric matrix, it makes sense to define maximal green sequences for such matrices.

\begin{dfn}

\rm {{Let $\widetilde{B}= \begin{pmatrix}
B\\I_{n} \end{pmatrix} \in M_{{2n}\times n}(\mathbb{Z})
$ be a {totally} sign-skew-symmetric matrix for which the {sign coherence property} holds and let $\widetilde{B}^{\sigma_s}= \begin{pmatrix}
B^{\sigma_s}\\
C^{\sigma_s}
\end{pmatrix}$ be the matrix obtained from $\widetilde{B}$ by a {composition of mutations $\mu_{\sigma_s}=\mu_{k_s}\mu_{k_{s-1}}....\mu_{k_1}$}, $1 \leq k_j \leq n$ for every $1 \leq j \leq s$,}  \begin{itemize}
                                                                                                                  \item {an index {$i$} in the matrix $\widetilde{B}^{\sigma_s}$ for $1\leq i \leq n$ is called  {\bf green }(respectively, {\bf red}) if the entries of the column indexed by $i$ in the $C-$matrix of $\widetilde{B}^{\sigma_s}$ are non-negative (respectively, non-positive).}
                                                                                                                  \item A sequence of indices {$(k_1,k_{2},....,k_{s})$}, where $ 1 \leq k_j \leq n$ for all $ j \in \{1,2,...,s\}$, is called a {\bf green sequence} if $k_{j}$ is green in the matrix $\mu_{k_{j-1}}....\mu_{k_{1}}(\widetilde{B})$ for $1 \leq j \leq s$. Such sequence is called { \bf maximal} if $\mu_{k_{s}}....\mu_{k_{1}}(\widetilde{B})$ does not have any green indices.

                                                                                                                \end{itemize}}

\end{dfn}

\begin{dfn}
\rm{A {\bf source} in a {sign-skew-symmetric} matrix $B$ of the size $n \times n$ is an index $i$ where $1 \leq i \leq n$ and {$b_{ik} \leq 0 $ }for all $1\leq k \leq n$ }. {In the associated simple quiver $\Delta(B)$ of $B$, the source $i$ has all the arrows incident to it emerging from it.}
\end{dfn}
\begin{dfn}
\rm{An {\bf admissible numbering by sources} of an acyclic sign-skew-symmetric matrix ${B}$ of the size $n \times n$ is an $n$-tuple $(i_1,i_2,...,i_n)$ such that the indices of $B$  are $\{i_1,...,i_n\}$ with $i_1$ a source in ${B}$ and the vertex $i_k$ is a source in $\mu_{i_{k-1}}....\mu_{i_1}({B})$ for any $2 \leq k \leq n$.}
\end{dfn}

\begin{lem}\label{source admissible}
\rm{Every {acyclic} sign-skew-symmetric matrix $B$ admits an admissible numbering by source.}
\begin{proof}
$B$ is a finite acyclic matrix, thus it has a source $i_1$. When mutating at this source, $\mu_{i_1}(b_{lj})=-b_{lj}$ if $l=i_1$ or $j=i_1$ and {$\mu_{i_1}(b_{lj})=b_{lj}$} otherwise. Let $\{B-{i_1}\}$ denote the matrix obtained from $B$ by deleting the $i_1$-th row and column, hence $\mu_{i_1}(\{B-{i_1}\})=\{B-{i_1}\}$. Since $B$ is acyclic, every submatrix is also acyclic. Then $\mu_{i_1}(\{B-{i_1}\})$ is also acyclic and thus it has a source $i_2\neq i_1$. {Again} the submatrix $\mu_{i_2}(\mu_{i_1}(\{\{B-{i_1}\}-{i_2}\}))$ {is the same submatrix {$\{B-i_1-i_2\}$} obtained from $B$ by deleting the rows and columns $i_1,i_2$}, {and it is acyclic} with a new source $i_3$. In every step the submatrix formed by {the indices} which haven't been mutated at is the same as the submatrix {obtained by the same {indices}} in the original one $B$ and the new source {$i_k$} in this submatrix is also a source in the whole matrix $\mu_{i_{k-1}}...\mu_{i_1}(B)$ for $1 \leq k \leq n-1$ since we are preforming mutations at sources so the {entries of the} other {indices} remain unchanged and {moreover the sign of} {entry $b_{i_ki_d}$} that connects the new source {$i_k$} with an old one {$i_d$} (which has already been mutated at) {{is non-positive in the matrix $\mu_{i_{k-1}}...\mu_{i_1}(B)$}} {where $1\leq d \leq k-1$}. By repeating this process $n-1$ times, the {index} $i_n$ will definitely be a source and we get an admissible numbering by source $(i_1,...,i_n)$.
\end{proof}
\end{lem}
$Q$ is a locally-finite quiver with a group of automorphisms $\Gamma$ and a {finite number of orbits $n$}.{{ If the index $i$ is a source in $Q$, then {$b_{ij}\leq 0$} for every $j \in Q_0$. Suppose \textcolor{blue}{$l \in \bar{i}$}, i.e there is an automorphism $g \in \Gamma$ such that {{$g(l)=i$}} and suppose that {$b_{l k}>0$} for some $k\in Q_0$, then by the definition of automorphism {$b_{l k}=b_{ig(k)} > 0 $}, contradiction}. Thus  any index in the orbit {$\bar{i}$ is also a source in $Q$. In this case we call $\bar{i}$ an {\bf{orbit-source}}} in $Q$.} If there is a sequence $(\bar{i}_1,....,\bar{i}_{n})$ of orbit-mutations such that the orbit $\bar{i}_1$ is an orbit-source in $Q$ and the orbit $\bar{i}_j$ is an orbit-source in $\mu_{\bar{i}_{j-1}}......\mu_{\bar{i}_{1}}(Q)$ for $1\leq j \leq n$ {and the set $\{\bar{i}_1,...,\bar{i}_n\}$ represents all the orbits under the action of $\Gamma$}, then the sequence $(\bar{i}_1,....,\bar{i}_{n})$ is called {\bf{orbit-admissible numbering by source}} in $Q$.\\
Getting back to the pair $(Q,\Gamma)$ as an unfolding of {the principal part} of an acyclic sign-skew-symmetric matrix $B$ constructed as in Construction \ref{cons2}. By the construction of $Q$, we notice that if a vertex labeled by $i$ is a source, then all the vertices labeled by $i$ are also sources.
\begin{cor}\label{source}
\rm{
Let $(Q,\Gamma)$ be the unfolding of {the principal part} of an acyclic sign-skew-symmetric matrix $B$ built as in Construction \ref{cons2} {with a finite set of orbits {$\{ \bar{1},..,\bar{n}\}$} obtained by the action of $\Gamma$}, and let $(\bar{i}_1,...,\bar{i}_k)$ be a sequence of orbit-mutation. If a vertex labeled by ${j}$ is a source in {$\mu_{\bar{i}_k}....\mu_{\bar{i}_1}(Q)$, then all the vertices labeled by $j$ are also sources in $\mu_{\bar{i}_k}....\mu_{\bar{i}_1}(Q)$ such that $\bar{i}_l \in \{\bar{1},...,\bar{n}\}$ for $1 \leq l \leq k$.}}
\begin{proof}
{Clearly $\Gamma$ is a group of automorphism for $\mu_{\bar{i}_k}....\mu_{\bar{i}_1}(Q)$. The statement holds true by the definition of automorphisms and since the indices which have the same label lie in the same orbit.}
\end{proof}
\end{cor}
\begin{lem}\label{orbitsource admissible}
\rm{Let $({Q},\Gamma)$ be the unfolding of {the principal part} of an acyclic sign-skew-symmetric matrix ${B} \in Mat_{n\times n}(\mathbb{Z})$ {built as in Construction \ref{cons2} with the adjacency matrix $B^\S$}, then ${Q}$ admits an orbit-admissible numbering by source.}
\begin{proof}
 By Lemma \ref{source admissible}, the matrix $B$ defines an admissible numbering by source $(i_1^f,...,{i_n}^f)$. {There are $n$ orbits obtained by the action of the group of automorphisms $\Gamma$ defined in Construction \ref{cons2}, each orbit has the vertices with the same label in $Q$ and since $\{i_1^f,...,{i_n}^f\}=\{1,2,3,...,n\}$, the set $\{\bar{i}_1,...,{\bar{i}_n}\}=\{\bar{1},\bar{2},\bar{3},...,\bar{n}\}$}. By Lemma \ref{invariantunfolding} {{the adjacency matrix of} the quiver $\mu_{\bar{i}_{j-1}}......\mu_{\bar{i}_1}(Q)$ is {an} unfolding of the matrix $\mu_{i_{j-1}^f}.....\mu_{i_1^f}(B)$ for $1 \leq j \leq n$. $i_j^f$ is a source in $\mu_{i_{j-1}^f}.....\mu_{i_1^f}(B)$, thus by the folding relation we get
  {$$\mu_{i_{j-1}^f}.....\mu_{i_1^f}(b_{i_j^fl^f})=\mu_{i_{j-1}^f}.....\mu_{i_1^f}(b_{\bar{i}_j\bar{l}})=\underset{r \in \bar{i}_j}{\sum}\mu_{\bar{i}_{j-1}}......\mu_{\bar{i}_1}(b_{rl}^\S)\leq 0$$}
 \\
  for every $ l \in Q_0$. {Since $\mu_{\bar{i}_{j-1}}......\mu_{\bar{i}_1}(Q)$ does not admit any $\Gamma-$2 cycles, the entries $\mu_{\bar{i}_{j-1}}......\mu_{\bar{i}_1}(b_{rl}^\S)$ have the same sign for every $r \in \bar{i}_j$ and every $l \in Q_0$}. Thus each term in the summation {above} is {non-positive} and hence $r$ is a source in $\mu_{\bar{i}_{j-1}}......\mu_{\bar{i}_1}(Q)$ for every $r \in \bar{i}_j$ and every $1 \leq j \leq n$. Therefore $(\bar{i}_1,...,\bar{i}_n)$ is an orbit-admissible numbering by source.}
 \end{proof}
 \end{lem}
{ By Lemma \ref{orbitlemmalocallyfinite}, we can define orbit-green sequences and orbit-maximal green green sequences for a locally-finite {framed} quiver $\widetilde{Q}$ with an adjacency matrix $\widetilde{B^\S}$ equipped with a group of automorphisms $\Gamma$ such that $\widetilde{Q}$ can do arbitrary steps of orbit-mutations.
 \begin{dfn}
  \rm{Let $\widetilde{Q}$ be a locally-finite {framed} quiver {with an adjacency matrix $\widetilde{B^\S}$ equipped with a group of automorphisms $\Gamma$ such that $\widetilde{Q}$ can do arbitrary steps of orbit-mutations}, and let $\bar{\sigma}_s :=(\bar{k}_1,\bar{k}_2,...,\bar{k}_s)$ be a sequence of orbit-mutations and let ${C^\S}^{\bar{{\sigma}}_{j}} $ be the $C$-$\textbf{matrix}$ of ${\widetilde{B^\S}}^{\bar{\sigma}_j}$ obtained from $\widetilde{B^\S}$ by the sequence of orbit-mutation $(\bar{k}_{1},\bar{k}_{2},.......,\bar{k}_{j})$ for $1 \leq j \leq s $, then $(\bar{k}_{1},\bar{k}_{2}.......,\bar{k}_{s})$ is said to be an {\bf orbit-green sequence} if for every $1\leq j \leq s $, $\bar{k}_{j}$ is a green orbit in ${B^\S}^{\bar{\sigma}_{j-1}}$. The sequence $(\bar{k}_{1},\bar{k}_{2}.......,\bar{k}_{s})$ is said to be {\bf orbit-maximal green sequence} if ${\widetilde{B^\S}}^{\bar{\sigma}_{s}}$ doesn't have any green orbits.}
\end{dfn}
}
{We always suppose that we have finitely many orbits under the action of $\Gamma$ on the unfolding locally-finite quiver of a sign-skew-symmetric matrix}.
\begin{lem}\label{superconnect}
\rm{Let {$\widetilde{Q}$ with the adjacency matrix} $\widetilde{B^\S}=\begin{pmatrix} B^\S\\ I_\infty \end{pmatrix}$ be the unfolding of an acyclic sign-skew-symmetric matrix $\widetilde{B}=\begin{pmatrix} B\\I_n
\end{pmatrix}$ {built as in Construction \ref{cons2}}, then any orbit admissible numbering by source {of} $\widetilde{B^\S}$ is an orbit-maximal green sequence}.
\begin{proof}
{Suppose that $\{\bar{i_1},.....,\bar{i_n}\}$ is an orbit-admissible numbering by source of $\widetilde{B^\S}$. By the definition of orbit-mutation, the mutation at a specific orbit-source reflects the arrows incident to the vertices of that orbit while keeping other arrows the same. Hence at each step we get a new red orbit while the colors of other orbits remain the same.}
\end{proof}
\end{lem}
{The following Theorem shows the relation between maximal green sequences for acyclic sign-skew-symmetric matrices and orbit-maximal green sequences for their unfolding matrices.}
\begin{mytheo}\label{connect}
\rm{{Let $\widetilde{B^\S}=\begin{pmatrix} B^\S\\ I_\infty \end{pmatrix}$ be the unfolding of an acyclic sign-skew symmetric matrix $\widetilde{B}=\begin{pmatrix} B\\I_n
\end{pmatrix}$ as constructed in Construction \ref{cons2}, then the sequence $(\bar{k}_{1},\bar{k}_{2}.......,\bar{k}_{s})$ is an orbit-maximal green sequence for $\widetilde{B^\S}$ if and only if the corresponding sequence $(k_{1}^f,k_{2}^f.......,k_{s}^f)$ is a maximal green sequence for its {folding} matrix $\widetilde{B}$.}}
\begin{proof}
{By Lemma \ref{invariantunfolding} the matrix $\mu_{\bar{k}_{j-1}}.....\mu_{\bar{k}_1}(\widetilde{B}^\S)$ is an unfolding of the matrix $\mu_{k_{j-1}^f}....\mu_{k_1^f}(\widetilde{B})$ for every $1\leq j \leq s$. Thus \begin{equation}\label{555}
\mu_{k_{j-1}^f}....\mu_{k_1^f}(c_{r^f{k}_j^f})=\underset{t\in \bar{r}}{\sum}\mu_{\bar{k}_{j-1}}.....\mu_{\bar{k}_1}(c_{tk_j}^\S)
\end{equation}
 The sign-coherence property is satisfied for locally-finite framed quivers by Lemma \ref{infinitesigncoherence}. Consequently, the terms that compose the summation on the right hand of equation (\ref{555}) have the same sign. Thus the $k_j^f$ is green (red) in the matrix $\mu_{k_{j-1}^f}....\mu_{k_1^f}(\widetilde{B})$ if and only if the index $k_j$ is green (red) in the matrix $\mu_{\bar{k}_{j-1}}.....\mu_{\bar{k}_1}(\widetilde{B}^\S)$ and equivalently by Lemma \ref{orbitlemmalocallyfinite} the orbit $\bar{k}_j$ is green (red) in the matrix $\mu_{\bar{k}_{j-1}}.....\mu_{\bar{k}_1}(\widetilde{B}^\S)$, hence the result follows.}
\end{proof}
\end{mytheo}
Now we can prove that every acyclic sign-skew-symmetric matrix admits a maximal green sequence.
\begin{mytheo}\label{mainn}
\rm{{Every acyclic sign-skew-symmetric matrix $\widetilde{B}=\begin{pmatrix} B\\I_n
\end{pmatrix}$ admits a maximal green sequence}.}
\begin{proof}
{Let $\widetilde{B^\S}=\begin{pmatrix} B^\S\\ I_\infty \end{pmatrix}$ be the unfolding of an acyclic sign-skew symmetric matrix $\widetilde{B}=\begin{pmatrix} B\\I_n
\end{pmatrix}$ as constructed in Construction \ref{cons2}.}
  By Lemma \ref{superconnect} the {locally-finite framed quiver $\widetilde{Q}$} with the {adjacency matrix} $\widetilde{B^\S}$ admits an orbit-maximal green sequence  $(\bar{i}_1,...,\bar{i}_n)$. Hence the sequence $({i}_1^f,...,{i}_n^f)$ is a maximal green sequence of matrix $\widetilde{B}$ by Theorem \ref{connect}.
\end{proof}
\end{mytheo}


\begin{thebibliography}{15}


\bibitem{AI}A. Assem, M. Blais, T. Brustle A. Samaon, \emph{Mutation Classes of Skew-Symmetric $3\times{3}$ Matrices}, Communications in Algebra, 36:4 (2008), 1209¨C1220.

\bibitem{nw} A. Berenstein, S. Fomin, A. Zelevinsky, \emph{Cluster algebras III: upper bound and Bruhat cells}, Duke
Math. J. 126:1 (2005), 1-52.

\bibitem{TB} T. Brustle, G. Dupton and M. Perotin, \emph{On Maximal Green Sequences},
Int. Math. Res. Not. 16 (2014),  4547-4586.

\bibitem{CAO}P. Cao and F. Li, \emph{ Uniformly column sign-coherence and the existence of maximal green sequences}, Journal of Algebraic Combinatorics 50 (2019),  403-417.

\bibitem {Derksen} H. Derksen, J. Weyman, A. Zelevinsky, \emph{Quivers with potentials and their representations
II: applications to cluster algebras}, J. Amer. Math. Soc. 23 (2010), 749-790.

\bibitem{dup}G. Dupton, \emph{An Approach to Non-Simply-Laced Cluster Algebra}, J. Algebra 320:4 (2008), 1626-1661.

\bibitem{FZ} S. Fomin and A. Zelevensky, \emph{Cluster Algebras. I.Foundations}, J. Amer. Math. Soc. 15:2 (2002),
	 497-529.

    \bibitem{SHAPIRO}S. Fomin, M. Shapiro and D. Thruston, \emph{Cluster Algrbra and Trangulated Surfaces Part I:Cluster Complex}, Acta Math. 201:1 (2008), 83-146.


\bibitem{GR}M. Gross, P. Hacking, S. Keel, M. Kontsevich, \emph{Canonical bases for cluster algebras},
J. Amer. Math. Soc. 31 (2018), 497-608.

\bibitem{mingli}M. Huang and F. Li,\emph{ Unfolding of sign-skew-symmetric cluster algebras
and its applications to positivity and F-polynomials}, Advances in Mathematics, 340 (2018),  221-283.

\bibitem{keller}B. Keller, \emph{Quiver Mutation and Combinatorial DT-Invariants}, corrected version of a contribution to FPSAC 2013, arXiv:1709.03143 [math.CO], 12 pages.

\bibitem{MILL}M. Mills, \emph{Maximal green sequences for quivers of finite mutation type}, Advances in Mathematics, 319 (2017), 182-210.

\bibitem{MULLER}G. Muller, \emph{The Existence of Maximal Green Sequence is Not Invariat Under Quiver Mutation}, The Electronic  Journal of Combinatorics,  23:2 (2016), Paper 2.47, 23pp.

\bibitem{NA}T. Nakanishi, S. Stella, \emph{Diagrammatic Description of c-Vectors and d-Vectors of Cluster Algebra of Finit Type}, The Electronic  Journal of Combinatorics,  21:1 (2014), Paper 1.3, 107pp.

\bibitem{SA}A. Seven, \emph{Maximal Green Sequences of Skew-Symmetrizable $3\times 3$ Matrices}, Linear Algebra and its Applications, 440 (2014),  125-130.



\end{thebibliography}
\end{document}